\documentclass[12pt, reqno]{amsart}
\usepackage{amsmath, amstext, amsbsy, amssymb, amscd}

\newtheorem{theorem}{Theorem}[section]
\newtheorem{lemma}[theorem]{Lemma}

\theoremstyle{definition}

\theoremstyle{remark}

\numberwithin{equation}{section}

\newcommand{\Hilb}{{\rm Hilb}}

\newcommand{\Z}{\mathbb{Z}}
\newcommand{\R}{\mathbb{R}}

\newcommand{\bZ}{\mathsf{Z}}

\newcommand{\Q}{\mathbb{Q}}

\begin{document}
\title[Donaldson-Thomas invariants]
{The Donaldson-Thomas invariants under blowups and flops}

\author[Jianxun Hu]{Jianxun Hu$^1$}
\address{Department of Mathematics, Zhongshan University,
Guangzhou, 510275, P. R. China} \email{stsjxhu@zsu.edu.cn}
\thanks{${}^1$Partially supported by the NSFC Grant 10231050 and NCET-04-0795}

\author[Wei-Ping Li]{Wei-Ping Li$^2$}
\address{Department of Mathematics, HKUST, Clear Water Bay,
Kowloon, Hong Kong} \email{mawpli@ust.hk}
\thanks{${}^2$Partially supported by the grant
HKUST6114/02P}

\begin{abstract}
Using the degeneration formula for Doanldson-Thomas invariants, we proved formulae
for blowing up a point and simple flops.

\end{abstract}

\maketitle
\date{}

\section{\bf Introduction}
Given a smooth projective Calabi-Yau $3$-fold $X$, the moduli space of stable
sheaves on $X$ has virtual dimension zero. Donaldson and Thomas \cite{D-T} defined
the holomorphic Casson invariant of $X$ which essentially counts the number of stable
bundles on $X$. However, the moduli space has positive dimension and is singular in
general. Making use of virtual cycle technique (see \cite{B-F} and \cite{L-T}),
Thomas showed in
\cite{Thomas} that one can define a  virtual moduli cycle for some $X$ including
Calabi-Yau and Fano $3$-folds. As a consequence, one can define Donaldson-type
invariants of $X$ which are deformation invariant. Donaldson-Thomas invariants
provide  a new  vehicle to study the geometry and other aspects of higher-dimensional
varieties. It is important to understand these invariants.

Much studied Gromov-Witten invariants of $X$ are the counting of stable maps from
curves to $X$. In \cite{MNOP1, MNOP2}, Maulik,  Nekrasov,  Okounkov, and
Pandharipande discovered relations between Gromov-Witten invarints of $X$ and
Donaldson-Thomas invariants constructed from moduli spaces of ideal sheaves of
curves on $X$. They conjectured that these two invariants can be identified via the
equations of partition functions of both theory. This suggests that many
phenomena on Gromov-Witten theory have the counterparts in Donaldson-Thomas theory.

Donaldson-Thomas invariants are deformation independent. In the birational geometry
of $3$-folds, we have blowups and flops. Donaldson-Thomas invariants couldn't be
effective in studying birational geometry unless we understand how invariants
change under birational operations. Li and Ruan in \cite{L-R} studied how
Gromov-Witten invariants change under a flop for Calabi-Yau $3$-fold. They proved
that one can identify the 3-point functions of $X$ and the flop $X^f$ of $X$ up to
some transformation of the $q$ variables. The same question was aslo studied by
Liu and Yau in \cite{L-Y} recently using the J. Li's degeneration formula from
algebraic geometry. In
\cite{Hu1, Hu2}, the first author studied the change of Gromov-Witten invariants
under the blowup. In this paper, we will study how Donaldson-Thomas invariants in
\cite{MNOP2} change under the blowup of a point and some flops.

The method we use is the degeneration formula for Donaldson-Thomas invariants
studied in \cite{Li1, Li2, MNOP2}.  The blowup of $X$ has a
description in terms of a degeneration of $X$. Then we can
apply the degeneration formula.  In the category of symplectic manifolds, one  uses
symplectic sum or symplectic cutting for the blowup operation on $X$.
The gluing formula for Gromov-Witten invariants in the
symplectic setup is in \cite{I-P1, I-P2, L-R}. Besides the difference of
degeneration and symplectic cutting, the arguments used in \cite{L-R, Hu1, Hu2, L-Y}
rely on the fact that stable maps have connected domain, while the curves defined
by ideal sheaves are in general not connected. Therefore the formula for the flop
is a bit different from that of Gromov-Witten invariants in
\cite{L-R}.

The organization of the paper is as follows. In section 2, we set up
terminologies and notations, and list the basic results needed. The degeneration
formula is discussed. In section 3, using J. Li's degeneration formula,  we
prove a blowup formula for the blowup of
$X$ at a point. In section 4, we prove the equality of Donaldson-Thomas partition
functions under a flop.

\bigskip\noindent
{\bf Acknowledgments.} Authors would like to thank Jun Li, Miles Reid, Qi Zhang,
Yongbin Ruan, and Zhenbo Qin for many helpful discussions. The second author would
like to thank the Department of Mathematics at Zhongshan University for the
hospitality during his several visits in the spring semester of 2005. The first
author would like to thank HKUST for the hospitality during his visit in
January of 2005. Both authors would like to thank the ICCM held at the Chinese
University of Hong Kong where they met and initiated the work.

\section{\bf Preliminaries}
In this section, we shall discuss the basic materials on Donaldson-Thomas
invariants studied by Maulik, Nekrasov, Okounkov and Pandharipande.
 For
the details, one can consult \cite{D-T, L-R, I-P1, I-P2, Li1, Li2, MNOP1, MNOP2,
Thomas}.

Let $X$ be a smooth projective  3-fold and $\mathcal I$ be an
ideal sheaf on $X$. Assume the sub-scheme $Y$ defined by $\mathcal
I$ has dimension $\le 1$. Here $Y$ is allowed to have embedded
points on the curve components. Therefore we have the exact
sequence
\begin{eqnarray*}
      0\longrightarrow {\mathcal I} \longrightarrow {\mathcal
      O}_X\longrightarrow {\mathcal O}_Y\longrightarrow 0.
\end{eqnarray*}

The $1$-dimensional components, with multiplicities taken into
consideration, determine a homology class
\begin{eqnarray*}
     [Y]\in H_2(X, \Z).
\end{eqnarray*}

Let $I_n(X,\beta)$ denote the moduli space of ideal sheaves
$\mathcal I$ satisfying
\begin{eqnarray*}
     \chi({\mathcal O}_Y) = n,
     \quad
   [Y] = \beta \in H_2(X, \Z).
\end{eqnarray*}
$I_n(X,\beta)$ is projective and is a  fine moduli space. From the
deformation theory, one can compute the virtual dimension of
$I_n(X,\beta)$ to obtain the following result

\begin{lemma}\label{lem2.1}  The virtual dimension of
$I_n(X,\beta)$, denoted by $\text{vdim}$,  equals $\int_\beta
c_1(T_X)$.
\end{lemma}

Note that the actual dimension of the moduli space $I_n(X,\beta)$
is usually larger than the virtual dimension.



 Let $\mathfrak I$ be the universal family
over $I_n(X,\beta)\times X$ and $\pi_i$ be the projection of
$I_n(X,\beta)\times X$ to the $i$-th factor. For a cohomology
class $\gamma \in H^l(X,\Z)$, consider the  operator
\begin{eqnarray*}
    ch_{k+2}(\gamma) : H_*(I_n(X,\beta), \Q) \longrightarrow
    H_{*-2k+2-l}(I_n(X,\beta),\Q),
\end{eqnarray*}
\begin{eqnarray*}
  ch_{k+2}(\gamma)(\xi) = \pi_{1*}(ch_{k+2}({\mathcal J})\cdot
  \pi_2^*(\gamma)\cap\pi_1^*(\xi)).
\end{eqnarray*}

Descendent fields in Donaldson-Thomas theory are defined in
\cite{MNOP2}, denoted by $\tilde{\tau}(\gamma)$, which correspond
to the operations $(-1)^{k+1}ch_{k+2}(\gamma)$. The descendent
invariants are defined  by
\begin{eqnarray*}
<\tilde{\tau}_{k_1}(\gamma_{l_1})\cdots
\tilde{\tau}_{k_r}(\gamma_{l_r})>_{n,\beta} =
\int_{[I_n(X,\beta)]^{vir}}\prod_{i=1}^r(-1)^{k_i+1}ch_{k_i+2}(\gamma_{l_i}),
\end{eqnarray*}
where the latter integral is the push-forward to a point of the
class
\begin{eqnarray*}
   (-1)^{k_1+1}ch_{k_1+2}(\gamma_{l_1})\circ \cdots \circ
   (-1)^{k_r+1}ch_{k_r+2}(\gamma_{l_r})([I_n(X,\beta)]^{vir}).
\end{eqnarray*}

The Donaldson-Thomas partition function with descendent insertions
is defined by
\begin{eqnarray*}
   Z_{DT}(X;q\mid\prod_{i=1}^r\tilde{\tau}_{k_i}(\gamma_{l_i}))_\beta
   = \sum_{n\in
   \Z}<\prod_{i=1}^r\tilde{\tau}_{k_i}(\gamma_{l_i})>_{n,\beta}q^n.
\end{eqnarray*}

The degree 0 moduli space $I_n(X,0)$ is isomorphic to the Hilbert
scheme of $n$ points on $X$. The degree 0 partition function is
$\bZ_{DT}(X;q)_0$.

The reduced partition function is obtained by formally removing
the degree $0$ contributions,
\begin{eqnarray*}
  Z'_{DT}(X;q\mid\prod_{i=1}^r\tilde{\tau}_{k_i}(\gamma_{l_i}))_\beta
  =
\frac{Z_{DT}(X;q\mid\prod\limits_{i=1}^r\tilde{\tau}_{k_i}(\gamma_{l_i}))_\beta}{Z_{DT}(X;q)_0}.
\end{eqnarray*}

Relative Donaldson-Thomas invarints are also defined in
\cite{MNOP2}. Let $S$ be a smooth divisor in $X$. An ideal sheaf
$\mathcal I$ is said to be relative to $S$ if the morphism
\begin{eqnarray*}
\mathcal I\otimes_{\mathcal O_X}\mathcal O_S\rightarrow \mathcal
O_X\otimes _{\mathcal O_X}\mathcal O_S
\end{eqnarray*}
is injective. A proper moduli space $I_n(X/S,\beta)$ of relative
ideal sheaves can be constructed by considering the ideal sheaves
relative to the expended pair $(X[k], S[k])$. For details, one can
read  \cite{Li2} and \cite{MNOP2}.

Let $Y$ be the subscheme defined by $\mathcal I$. The scheme
theoretic intersection $Y\cap S$ is an element in the Hilbert
scheme of points on $S$ with length $[Y]\cdot S$. If we use
$\Hilb(S, k)$ to denote the Hilbert scheme of points of
length $k$ on $S$, we  have a map
\begin{eqnarray*}
   \epsilon : I_n(X/S,\beta) \longrightarrow \Hilb(S,
   \beta\cdot [S]).
\end{eqnarray*}

The cohomology of the Hilbert scheme of points of $S$ has a basis
via the representation of the Heisenberg algebra on the
cohomologies of the Hilbert schemes.

Following Nakajima in \cite{Nakajima}, let $\eta$ be a cohomology weighted partition
with respect to a basis of $H^*(S, \Q)$. Let $\eta=\{\eta_1, \ldots, \eta_s\}$ be a
partition whose corresponding cohomology classes are $\delta_1,
\cdots, \delta_s$, let
\begin{eqnarray*}
    C_\eta
    =\frac{1}{\mathfrak{z}(\eta)}P_{\delta_1}[\eta_1]\cdots P_{\delta_s}[\eta_s]\cdot{\bf
    1}\in H^*(\mbox{Hilb}(S, |\eta|), \Q),
\end{eqnarray*}
where
\begin{eqnarray*}
     \mathfrak{z}(\eta) = \prod_i \eta_i|\mbox{Aut}(\eta)|,
\end{eqnarray*}
and $|\eta|=\sum_j\eta_j$. The Nakajima basis of the cohomology of
$\mbox{Hilb}(S,k)$ is the set,
\begin{eqnarray*}
   \{C_\eta\}_{|\eta|=k}.
\end{eqnarray*}

We can choose a basis of $H^*(S)$ so that it is self dual with
respect to the Poincar\'e pairing, i.e., for any $i$, $\delta_i^*
=\delta_j$ for some $j$.  To each weighted partition $\eta$, we
define the dual partition $\eta^\vee$ such that
$\eta^\vee_i=\eta_i$ and the corresponding cohomology class to
$\eta^\vee_i$ is $\delta_i^*$. Then we have
\begin{eqnarray*}
 \int_{\mbox{Hilb}(S,k)}C_\eta\cup C_\nu =
 \frac{(-1)^{k-\ell(\eta)}}{\mathfrak{z}(\eta)}\delta_{\nu,\eta^\vee},
\end{eqnarray*}
see \cite{Nakajima}.

The descendent invariants in the relative Donaldson-Thomas theory are
defined by
\begin{eqnarray*}
<\tilde{\tau}_{k_1}(\gamma_{l_1})\cdots\tilde{\tau}_{k_r}(\gamma_{l_r})\mid
\eta>_{n,\beta} =
\int_{[I_n(X/S,\beta)]^{vir}}\prod_{i=1}^r(-1)^{k_i+1}ch_{k_i+2}(\gamma_{l_i})\cap
\epsilon^*(C_\eta),
\end{eqnarray*}

Define the associated partition function by
\begin{eqnarray*}
   Z_{DT}(X/S;q\mid\prod_{i=1}^r\tilde{\tau}_{k_i}(\gamma_{l_i}))_{\beta,\eta}
   = \sum_{n\in\Z}
   <\prod_{i=1}^r\tilde{\tau}_{k_i}(\gamma_{l_i})\mid\eta>_{n,\beta}q^n.
\end{eqnarray*}

The reduced partition function is obtained by formally removing
the degree $0$ contributions,
\begin{eqnarray*}
Z'_{DT}(X/S;q\mid\prod_{i=1}^r\tilde{\tau}_{k_i}(\gamma_{l_i}))_{\beta,\eta}=
\frac{Z_{DT}(X/S;q\mid\prod\limits_{i=1}^r\tilde{\tau}_{k_i}(\gamma_{l_i}))_{\beta,\eta}}{Z_{DT}(X/S;q)_0}.
\end{eqnarray*}

In the remaining of the section, we shall discuss the degeneration formula due to
J. Li. It is the main tool employed in the paper.

 Let $\pi\colon\mathcal X\to C$ be
a smooth
$4$-fold over a smooth irrreducible curve $C$ with a marked point denoted by $\bf 0$
such that $\mathcal X_t=\pi^{-1}(t)\cong X$ for $t\neq {\bf 0}$ and
$\mathcal X_{\bf 0}$ is a union of two smooth $3$-folds $X_1$ and
$X_2$ intersecting transversely along a smooth surface $S$. We
write $\mathcal X_{\bf 0}=X_1\cup_SX_2$. Assume that $C$ is
contractible and $S$ is simply-connected.

Consider the natural maps
\begin{eqnarray*}
i_t\colon X=\mathcal X_t\rightarrow \mathcal X,\qquad i_{\bf
0}\colon \mathcal X_{\bf 0}\rightarrow \mathcal X,
\end{eqnarray*}
and the gluing map
\begin{eqnarray*}
g=(j_1, j_2)\colon X_1\coprod X_2\rightarrow \mathcal X_{\bf 0}.
\end{eqnarray*}

We have
\begin{eqnarray*}
H_2(X){\buildrel{i_{t*}}\over \longrightarrow} H_2(\mathcal
X){\buildrel{i_{0*}}\over \longleftarrow} H_2(\mathcal X_{\bf
0}){\buildrel{g_*}\over\longleftarrow} H_2(X_1)\oplus H_2(X_2),
\end{eqnarray*}
where $i_{0*}$ is an isomorphism since there exists a deformation
retract from $\mathcal X$ to $\mathcal X_{\bf 0}$ (see \cite{Clemens}) and $g_*$ is
surjective from Mayer-Vietoris sequence. For $\beta \in H_2(X)$,
there exist $\beta_1\in H_2(X_1)$ and $\beta_2\in H_2(X_2)$ such
that
\begin{eqnarray}\label{betasum}
i_{t*}(\beta)=i_{0*}(j_{1*}(\beta_1)+j_{2*}(\beta_2)).
\end{eqnarray}
For simplicity, we write $\beta=\beta_1+\beta_2$ instead.
\begin{lemma}\label{lem2.3}
With the assumption as above,  given $\beta = \beta_1 + \beta_2$.  Let $
d = \int_\beta c_1(X)$ and $d_i = \int_{\beta_i}c_1(X_i)$, $ i =
1,2$. Then
\begin{eqnarray}\label{deg-formula}
    d = d_1 + d_2 - 2\int_{\beta_1}[S],\qquad  \int_{\beta_1} [S]=\int_{\beta_2}
    [S].
\end{eqnarray}
\end{lemma}
\begin{proof}
The formulae (\ref{deg-formula}) come from the adjunction formulae
$K_{\mathcal X_t}=K_{\mathcal X}|_{\mathcal X_t}$ and
$K_{X_i}=(K_{\mathcal X}+X_i)|_{X_i}$ for $i=1, 2$, and $X_1\cdot
(X_1+X_2)=X_1\cdot \mathcal X_{\bf 0}=0$.
\end{proof}


 Similarly for cohomology, we have the maps
\begin{eqnarray*}
H^k(\mathcal X_t){\buildrel{i_t^*}\over\longleftarrow}
H^k(\mathcal X){\buildrel{i_0^*}\over \longrightarrow }
H^k(\mathcal X_{\bf
0}){\buildrel{g^*}\over\longrightarrow}H^k(X_1)\oplus H^k(X_2),
\end{eqnarray*}
where $i_0^*$ is an isomorphism. Take $\alpha\in H^k(\mathcal X)$
and let
$\alpha(t)=i^*_t\alpha$.

There is a degeneration formula  which takes the form
\begin{eqnarray}\label{gluing-formula}
 &Z'_{DT}(\mathcal
X_t;q\mid\prod\limits_{i=1}^r\tilde{\tau}_0(\gamma_{l_i}(t)))_\beta\nonumber\\
 =
&\sum
Z'_{DT}({X_1}/{S};q\mid\prod\tilde{\tau}_0(j_1^*\gamma_{l_i}(0)))_{\beta_1,\eta}\displaystyle{
\frac{(-1)^{|\eta|-\ell(\eta)}\mathfrak{z}(\eta)}{q^{|\eta|}}}
\nonumber\\
&\cdot
Z'_{DT}({X_2}/{S};q\mid\prod
\tilde{\tau}_0(j_2^*\gamma_{\l_i}(0)))_{\beta_2,\eta^\vee},
\end{eqnarray}
where the sum is over the splittings $\beta_1 + \beta_2 = \beta$,
and cohomology weighted partitions $\eta$. $\gamma_{l_i}$'s are
cohomology classes on $\mathcal X$. There is a compatibility
condition
\begin{eqnarray}\label{comp}
|\eta|=\beta_1\cdot [S]=\beta_2\cdot [S].
\end{eqnarray}

For details, one can see \cite{Li1, Li2, MNOP2}.

\section{\bf Blowup at a point and a Blowup formula}

   In \cite{MNOP1, MNOP2}, the authors discovered a correspondence
between Gromov-Witten theories and Donaldson-Thomas theories. In
\cite{Hu1, Hu2}, the first author studied the change of
Gromov-Witten invariants under the blowup operation. In this
section, we will study the change of Donaldson-Thomas invariants
under the blowup along a point.

    The key idea  is that the blowup can be obtained via a semistable
degeneration as follows. Let $X$ be a smooth projective
 3-fold and $\tilde{X}$ be the blowup of $X$ at a general
point $x$. Denote by $p:\tilde{X}\longrightarrow X$ the natural
projection of the blowup.   Let $\mathcal X$ be the blow up of
$X\times \mathbb C$ at the point $(x, 0)$ and let $\pi$ be the
natural projection from $\mathcal X$ to $\mathbb C$. It is a
semistable degeneration of $X$ with the central fiber $\mathcal
X_0$ being a union of $X_1 \cong \tilde{X}$ and $X_2 \cong \mathbb
P^3$, which is the exceptional divisor in $\mathcal X$. $X_1$ and
$X_2$ intersect transversely along $E\cong \mathbb P^2$, which is
the exceptional divisor in $X_1=\tilde X$. As a divisor in $X_2$, $E$ is a
hyperplane. $c_1(X_2)=4E$.

\begin{theorem}\label{blowupthm}
 Let $X$ be a smooth projective
3-fold.  Suppose that $\beta \in H_2(X, \Z)$ and
$\gamma_{l_i}\in H^*(X, \R)$, $i = 1, \cdots, r$. Then
\begin{eqnarray}\label{blowupformula}
  Z'_{DT}(X;q\mid\prod_{i=1}^r\tilde{\tau}_{0}(\gamma_{l_i}))_{\beta} =
Z'_{DT}(\tilde{X};q\mid\prod_{i=1}^r\tilde{\tau}_{0}(p^*\gamma_{l_i}))_{p^!(\beta)},
\end{eqnarray}
where  $p^!(\beta) =
PDp^*PD^{-1}(\beta)$.
\end{theorem}
\begin{proof}
 Choose  the support of
$\gamma_{l_i}$ outside of $x$. Then we have $\gamma_{l_i}\in
H^*(X_1)$ and no $\gamma_{\l_i}$'s in $H^*(X_2)$. In fact, let $p_1\colon \mathcal
X\to X$ be the composition of the blowing-down map $\mathcal X\to X\times \mathbb C$
with the projection $X\times\mathbb C\to X$. One can check that
$i^*_tp_1^*\gamma_{l_i}=\gamma_{l_i}$ and
$j_1^*i_0^*p_1^*\gamma_{l_i}=p^*\gamma_{l_i}$ and
$j_2^*i_0^*p_1^*\gamma_{l_i}=0$. We apply the degeneration formula
(\ref{gluing-formula}) to the cohomology classes $p_1^*\gamma_{\ell_i}$ on
$\mathcal X$.

   By the degeneration formula (\ref{gluing-formula}), we may express the absolute
Donaldson-Thomas invariants of $X$ in term of the relative
Donaldson-Thomas invariants of $(X_1, E)$ and $(X_2,E)$ as
follows:
\begin{eqnarray}\label{form1}
&Z'_{DT}(X;q\mid\prod\limits_{i=1}^r\tilde{\tau}_{0}(\gamma_{l_i}))_{\beta}\\
   = &   \sum\limits_{\eta, \beta_1 + \beta_2=
  \beta}Z'_{DT}\big({X_1}/{E};q\mid
  \prod\limits_{i=1}^r\tilde{\tau}_0(p^*\gamma_{l_i})\big)_{\beta_1,\eta}
\displaystyle{\frac{(-1)^{|\eta|-\ell(\eta)}\mathfrak{z}(\eta)}{q^{|\eta|}}}
Z'_{DT}({X_2}/{E};q)_{\beta_2,\eta^\vee}.\nonumber
\end{eqnarray}

 Now we need to compute the summands in the right hand side of the
degeneration formula. For this we have the following claim:

\noindent {\bf Claim:} There are only terms with $\beta_2 = 0$.

In fact, if $|\eta|\not=0$, then $\beta_2 \not= 0$ because
$\beta_2\cdot E =|\eta|$. By Lemma \ref{lem2.1},  we have
\begin{eqnarray*}
  c_1(X_1)\cdot \beta_1
= \mbox{vdim} I_n({X_1}/{E},\beta_1) = \sum_{i=1}^r \deg ch_2(\gamma_{l_i}) + \deg
  \epsilon_1^*(C_\eta),
\end{eqnarray*}
where $\epsilon_1 : I_n({X_1}/{E},\beta_1) \longrightarrow
\mbox{Hilb}(E,|\eta|)$ is the canonical intersection map, and
\begin{eqnarray*}
c_1(X_2)\cdot \beta_2 &=&\mbox{vdim} I_n({X_2}/{E},\beta_1)  =
4E\cdot \beta_2 = 4|\eta|,\\
   c_1(X)\cdot \beta& =& \mbox{vdim}
I_n(X,\beta) = \sum\limits_{i=1}^r \deg ch_2(\gamma_{l_i}).
\end{eqnarray*}
We have the last equality above because, otherwise, the involved
Donaldson-Thomas invariants of $X$ and $\tilde{X}$ will vanish and
the theorem holds.

By (\ref{deg-formula}), we have
\begin{eqnarray*}
  c_1(X)\cdot \beta = c_1(X_1)\cdot \beta_1 + c_1(X_2)\cdot
  \beta_2 - 2|\eta|.
\end{eqnarray*}
Combining all the four equations above, we obtain
\begin{eqnarray*}
  0 = \deg C_\eta + 2 |\eta|.
\end{eqnarray*}
This is a contradiction. Therefore $|\eta|=0$. So the claim is proved.

Thus $\beta_2\cdot E=0$. Since $E$ is the hyperplane in $X_2\cong \mathbb P^3$, we
must have $\beta_2=0$.   Also we have $\beta_1=p^!(\beta)$.

By the degeneration formula, we have
\begin{eqnarray}\label{blowup1}
  & &Z'_{DT}(X;q\mid
  \prod_{i=1}^r\tilde{\tau}_0(\gamma_{l_i}))_\beta\nonumber\\
=&&Z'_{DT}({X_1}/{E};q\mid
\prod_{i=1}^r\tilde{\tau}_0(p^*\gamma_{l_i}))_{p^!(\beta)}.
\end{eqnarray}

Now we want to use the degeneration formula one more time to study
the Donaldson-Thomas invariants of $\tilde{X}$. We blow up $\tilde
X\times \mathbb C$ along the surface $E\times 0$ to get a $4$-fold
$\tilde{\mathcal X}$. There is a projection $\tilde{\pi}\colon
\tilde{\mathcal X}\to \mathbb C$. The central fiber is a union of
$\tilde{X_1}=\tilde X$ and $\tilde{X_2}=\mathbb P(\mathcal
O_E(-1)\oplus \mathcal O_E)$ intersecting transversely along a
smooth surface $Z$, which is the surface $E$ in $\tilde{X_1}$ and
the infinite section $D_{\infty}$ in the projective bundle
$\tilde{X_2}$. Note that $\tilde{X_2}-D_{\infty}$ is the line
bundle $\mathcal O_E(-1)$, $p^!(\beta)\cdot E=0$, and
$PD(\gamma_{l_i})\cap E=\emptyset$. Let $\tilde{p_1}$ be the
composition of the map $\tilde{\mathcal X}\to \tilde X\times
\mathbb C$ and the map $\tilde X\times \mathbb C\to \tilde X$.
Applying the degeneration formula (\ref{gluing-formula}) to the
cohomology classes $\tilde{p_1}^*(\gamma_{l_i})$, we have
\begin{eqnarray*}
  &&Z'_{DT}(\tilde{X};q\mid
  \prod\limits_{i=1}^r\tilde{\tau}_0(\gamma_{l_i}))_{p^!(\beta)}\\
  =&&
 \sum_{\beta_1 + \beta_2 =
 p^!(\beta),\,\eta}Z'_{DT}({\tilde{X}_1}/{Z};q\mid
 \prod\limits_{i=1}^r\tilde{\tau}_0(p^*\gamma_{l_i}))_{\beta_1,
 \eta}\displaystyle{\frac{(-1)^{|\eta|-\ell(\eta)}\mathfrak{z}(\eta)}{q^{|\eta|}}}
 Z'_{DT}({\tilde{X}_2}/{Z};q)_{\beta_2, \eta^\vee},\nonumber
 \end{eqnarray*}
where $\beta_1 \cdot Z = |\eta|$.

Here we have the following claim as in the first part of our
proof:

\noindent{\bf Claim:} There are only terms with $\beta_2 = 0$ and
no $\eta$.

It is easy to see that $\tilde{X}_2
$ is the blowup $\tilde{\mathbb P^3} $ of ${\mathbb P}^3$ at a
point $p_0$. Denote by $\rho : \tilde{{\mathbb P}^3}
\longrightarrow {\mathbb P}^3$ the projection of the blowup. Let
$\ell\subset \tilde{\mathbb P^3}$ be the strict transform of a line in ${\mathbb
P}^3$ passing through the blown-up point $p_0$, and $S$ be the exceptional surface
of the blowup $\rho$. Denote by $e$ a line in $S$ which is an extremal ray. Since
$\ell$ is a fiber of ${\mathbb P}(\mathcal O_E(-1)\oplus \mathcal O_E)
\longrightarrow E$ which  also is an extremal ray, by Mori's theory, we have
$\beta_2 = a \ell + b e$,
$ a\geq 0$, $b\geq 0$. Let $H$ be the hyperplane class in $\mathbb P^3$. Since
$\rho^*H
\sim D_\infty$,  we have
$ a = \rho^*H\cdot \beta_2 =|\eta| $. One can show that $2D_{\infty}\cdot
\beta_2=p^!(\beta)\cdot E$. Since
$p^!(\beta)\cdot E = 0$, we have $a=|\eta|=\beta_2\cdot
D_{\infty}=0$ and $\beta_1\cdot E=\beta_2\cdot
D_{\infty}=0$.

We have
\begin{eqnarray*}
p^!(\beta)=\tilde{p_1}_{*}(p^!\beta)=\tilde{p_1}_*(\beta_1+\beta_2)=\beta_1+be,
\end{eqnarray*}
where we still use the same $e$ to represent a  line in $E\subset \tilde X$.

Since $E\cdot p^!(\beta)=0$ and $E\cdot \beta_1=0$, we have $-b=be\cdot E=0$. Thus
$b=0$ and hence $\beta_2=0$. The claim is proved.

We also see that $\beta_1=p^!(\beta)$.






By  the degeneration formula, we have
\begin{eqnarray}\label{blowup2}
  Z'_{DT}(\tilde{X};q\mid
  \prod\limits_{i=1}^r\tilde{\tau}_0(p^*\gamma_{l_i}))_{p^!(\beta)} & = &
 Z'_{DT}({\tilde{X}_1}/{Z};q\mid
 \prod\limits_{i=1}^r\tilde{\tau}_0(p^*\gamma_{l_i}))_{p^!(\beta)}\cdot
  Z'_{DT}({\tilde{X}_2}/{Z};q)_0\nonumber\\
  & = & Z'_{DT}({\tilde{X}}/{E};q\mid
 \prod_{i=1}^r\tilde{\tau}_0(p^*\gamma_{l_i}))_{p^!(\beta)}.
\end{eqnarray}

Note that $\tilde X_1\cong \tilde X$. Comparing (\ref{blowup1}) with
(\ref{blowup2}), we proved the Theorem.
\end{proof}


\section{\bf Blowup of $(-1, -1)$-curves and a flop formula}

In this section, we will study how Donaldson-Thomas invariants change under some
flops. The materials related to the birational geometry of $3$-folds can be found in
\cite{Kollar}, \cite{Kawamata}, \cite{KMM}, \cite{K-M}, \cite{Matsuki}.

 Let
$X$ be a smooth projective
$3$-fold,
$D$ be an effective divisor on
$X$. Suppose that $X$ admits a contractoin
of an extremal ray with respect to
$K_X+\epsilon D$, where $0<\epsilon \ll 1$,
\begin{eqnarray*}
\varphi\colon X\longrightarrow Y.
\end{eqnarray*}

Assume furthermore that the exceptional locus $Exc(\varphi)$ of $\varphi$ consists
of finitely many disjoint smooth rational $(-1, -1)$-curves $\Gamma_2, \ldots,
\Gamma_{\ell}$. $Y$ is a normal projective variety, $-D$ is $\varphi$-ample, and all
curves $\Gamma_i$ are numerically equivalent. Let's use $[\gamma]$ to denote the
numerically equivalent classes $\Gamma_i$, $i=2, \ldots, \ell$. There exists a
smooth projective
$3$-fold $X^f$ and a morphism
\begin{eqnarray*}
\varphi^f\colon X^f\longrightarrow Y,
\end{eqnarray*}
which is the flop of $\varphi$. $X^f$ can be obtained as follows in our situation.
We blow up $X$ along all the curves $\Gamma_i$, $i=2, \ldots, \ell$ to get a smooth
projective $3$-fold $\widetilde X$ with the exceptional divisors
$E_i\cong\Gamma_i\times \mathbb P^1$, $i=2, \ldots, \ell$. Let $\mu\colon \widetilde
X\rightarrow X$ be the blowup map.  We can blow down
$\widetilde X$ along all the $\Gamma_i$-direction. The new $3$-fold $X^f$ is
smooth, projective and containing $(-1, -1)$-curves $\Gamma_i^f$ for $i=2, \ldots,
\ell$. $\Gamma_i^f$ is the image of $E_i$ under the blow down. $X$ and $X^f$ are
birational and isomorphic in codimension one.

For any divisor $B$ on $X$, let $B^f$ be the strict transform of $B$ in $X^f$. We
have an isomorphism $N^1(X)\cong N^1(X^f)$ and
\begin{eqnarray*}
N^1(X)\cong \varphi^*N^1(Y)\oplus \mathbb R[D],\quad
N^1(X^f)\cong (\varphi^f)^*N^1(Y)\oplus \mathbb R[D^f].
\end{eqnarray*}

Similarly we get an isomorphism $H_2(X)\rightarrow H_2(X^f)$, denoted by $\phi_*$,
such that $\phi_*([\Gamma_i])=-[\Gamma_i^f]$ (see \cite{L-R}). The map $\phi_*$
induces isomorphisms
$\phi^*\colon H^{2i}(X^f)\rightarrow H^{2i}(X)$.

The map $\phi_*$ can also be seen as follows (see \cite{L-R}). There is an
injection $\iota$ from $H_2(X)$ to $H_2(\tilde X)$ such that the image of $\iota$
is the set $\{ \beta\in H_2(\tilde X)\,|\, \beta\cdot E=0\}$ where $E$ is the
exceptional divisor of the blow up. Similarly, there is an injection $\iota^f$ from
$H_2(X^f)$ to $H_2(\tilde X)$ with the same image. In fact, $(\iota^f)^{-1}\circ
\iota$ induces the isomorphism $\phi_*$.


 Let
$\mathcal X$ be the blow up of
$X\times
\mathbb C$ along all the curves
$\Gamma_i\times 0$. Let $\pi\colon \mathcal X\to \mathbb C$ be the natural
projection. Thus we get a semi-stable degeneration of $X$ whose central
fiber is a union of $X_1\cong \tilde X$ and $X_i=\mathbb P(\mathcal
O_{\Gamma_i}(-1)\oplus\mathcal O_{\Gamma_i}(-1)\oplus \mathcal O_{\Gamma_i})$ for
$i=2, \ldots, \ell$ with $X_1$ and $X_i$ intersecting transversely along the smooth
surface
$E_i$.

Here is a technical lemma.
\begin{lemma}\label{analytic-continuation}
The power series $\sum\limits_{d>0}d^kx^d$ has an analytic continuation $f_k(x)$ in
the domain
$\mathbb C-\{1\}$ such that
\begin{eqnarray*}
f_k({x}^{-1})=(-1)^{k+1}f_k(x).
\end{eqnarray*}
\end{lemma}
\begin{proof}
>From the geometric series formula $1+x+\ldots+x^d+\ldots=(1-x)^{-1}$, we get
\begin{eqnarray*}
x+2x^2+\ldots+dx^d+\ldots=x\cdot (1+x+\ldots+x^d+\ldots)^{\prime}=\frac{x}{(1-x)^2}.
\end{eqnarray*}
Let $f_1(x)=\displaystyle\frac{x}{(1-x)^2}$. One can check that
$f_1(x^{-1})=f_1(x)$.

Assume that the statement in the Lemma holds for $k$. Then
\begin{eqnarray*}
x+2^{k+1}x^2+\ldots+d^{k+1}x^d+\ldots=x\cdot
(x+\ldots+d^kx^d+\ldots)^{\prime}
\end{eqnarray*}
has an analytic continuation $f_{k+1}(x)=f^{\prime}_k(x)\cdot x$.
>From the chain rule, one has
$f_k^{\prime}(x^{-1})(-x^{-2})=(-1)^{k+1}f^{\prime}_k(x)$. Therefore
\begin{eqnarray*}
f_{k+1}(x^{-1})=x^{-1}f^{\prime}_{k}(x^{-1})=(-1)^{k+2}xf^{\prime}_k(x)=(-1)^{k+2}f_{k+1}(x).
\end{eqnarray*}
By the mathematical induction, we proved the Lemma.
\end{proof}

>From the proof, one can see that $f_k(x)=f_k(x^{-1})$ when $k$ is odd.

Define a function $g(q, v, \Gamma)$ by
\begin{eqnarray}\label{function-g}
g(q, v, \Gamma)=exp\{u^{-2}\sum_{d>0}\frac{1}{d^3}v^{d\Gamma}\}\cdot
\frac{1}{(1-v^{\Gamma})^{1/12}},
\end{eqnarray}
where $q=-e^{iu}$.

\begin{theorem}\label{flopthm}
Suppose cohomology classes
$\gamma_{l_i}\in H^{2k}(X^f)$,
$i=1,\cdots, r$ and $k=1,2,3$, have supports away from all the exceptional
curve
$\Gamma_i$.
\begin{enumerate}
\item[(i)] If $\beta = m[\gamma]$, we have
\begin{eqnarray*}
  Z'_{DT}(X;q)_{\beta} =
  Z'_{DT}(X^f;q)_{-\phi_*(\beta)}.
\end{eqnarray*}
\item[(ii)] There exist power series
\begin{eqnarray*}\Phi_X(q,
v|\{\phi^*\gamma_{\ell_i}\})&=&\sum\limits_{\beta\in
\iota(H_2(X))}\Phi_X(q |\{\phi^*\gamma_{\ell_i}\})_{\beta}\cdot v^{\beta},\\
\Phi_{X^f}(q,
v|\{\gamma_{\ell_i}\})&=&\sum\limits_{\beta\in
\iota^f(H_2(X^f))}\Phi_{X^f}(q |\{\gamma_{\ell_i}\})_{\beta}\cdot v^{\beta},
\end{eqnarray*}
 and
$G(q, v,\Gamma)$ such that
\begin{eqnarray*}
\Phi_X(q,
v|\{\phi^*\gamma_{\ell_i}\})
=
\Phi_{X^f}(q,
v|\{\gamma_{\ell_i}\})
,
\end{eqnarray*}
$G(q, v, \Gamma_i)/g(q, v, \Gamma_i)$ and $G(q, v^{-1}, \Gamma^f_i)/g(q, v^{-1},
\Gamma^f_i)$ are equivalent under analytic continuation, and
\begin{eqnarray}\label{flop-split1}
  Z'_{DT}(X;q,v\mid
\prod_{i=1}^r\tilde{\tau_0}(\phi^*\gamma_{l_i}))
=\Phi_X(q,
v|\{\phi^*\gamma_{\ell_i}\})\cdot\prod_{i=1}^{\ell}G(q,v,\Gamma_i),
\end{eqnarray}

\begin{eqnarray}\label{flop-split2}
  Z'_{DT}(X^f;q,v\mid
\prod_{i=1}^r\tilde{\tau_0}(\gamma_{l_i}))
=\Phi_{X^f}(q,
v|\{\gamma_{\ell_i}\})\cdot\prod_{i=1}^{\ell}G(q,v,\Gamma^f_i).
\end{eqnarray}
\end{enumerate}
\end{theorem}

\begin{proof}
There is a degeneration formula similar to (\ref{gluing-formula}) (see
\cite{Li2}) for the degeneration $\mathcal X$ described above. For simplicity, we
shall prove the case when there is only one $\Gamma_i$, denoted by $\Gamma$. The
proof for the general case is similar.

 By the degeneration formula (\ref{gluing-formula}),
we have
\begin{eqnarray*}
  &&Z'_{DT}\big(X;q\mid
  \prod\limits_{i=1}^r\tilde{\tau_0}(\phi^*\gamma_{l_i})\big)_\beta\\
 =&&
  \sum\limits_{\eta, \beta_1 + \beta_2 = \beta}Z'_{DT}\big({X_1}/{E};q\mid
\prod\limits_{i=1}^r\tilde{\tau_0}(\mu^*\phi^*\gamma_{l_i})\big)_{\beta_1,\eta}
\displaystyle{\frac{(-1)^{|\eta|-\ell(\eta)}\mathfrak{z}}
  {q^{|\eta|}}}Z'_{DT}({X_2}/{E};q)_{\beta_2, \eta^\vee},
\end{eqnarray*}
where $E$ is the intersection of $X_1$ with $X_2$, which is also the exceptional
divisor in $X_1$.

Similar to the proof of Theorem \ref{blowupthm}, we need to study the summands in
RHS. Therefore, we also need to compute the virtual dimensions of
involved moduli spaces. About the contributions of each term in
RHS, we have the following claim:

\noindent{\bf Claim:} There are only terms without
$\eta$.

In fact, suppose that $|\eta|\not = 0$. First of all, we want to
compute the first Chern class of $X_2$.

Let $V = {\mathcal O}_\Gamma(-1)\oplus {\mathcal O}_\Gamma(-1)\oplus
{\mathcal O}_\Gamma$ and  $p: \mathbb P(V)\longrightarrow \Gamma$
be the  projection. $X_2=\mathbb P(V)$.
For this projective bundle, we have the Euler exact sequence
\begin{eqnarray*}
  0 \longrightarrow {\mathcal O}_{{\mathbb P}(V)}\longrightarrow
  p^*V\otimes {\mathcal O}_{{\mathbb P}(V)}(1) \longrightarrow T_{{\mathbb
  P}(V)/\Gamma} \longrightarrow 0,
\end{eqnarray*}
 We also have
\begin{eqnarray*}
 0 \longrightarrow p^*\Omega_\Gamma^1 \longrightarrow \Omega_{{\mathbb
 P}(V)}^1\longrightarrow \Omega_{{\mathbb P}(V)/\Gamma}^1 \longrightarrow
 0.
\end{eqnarray*}
Therefore, we have
\begin{eqnarray*}
c_1(\Omega_{{\mathbb P}(V)}^1) & = & p^*c_1(\Omega_\Gamma^1) +
c_1(\Omega_{{\mathbb P}(V)/\Gamma}^1)\\
& = & p^*c_1(\Omega_\Gamma^1) - c_1(p^*V\otimes {\mathcal O}_{{\mathbb
P}(V)}(1))\\
& = & p^*c_1(K_\Gamma) - p^*c_1(V) - 3 c_1({\mathcal O}_{{\mathbb
P}(V)}(1))\\
& = & -3 c_1({\mathcal O}_{{\mathbb P}(V)}(1)),
\end{eqnarray*}
where $c_1({\mathcal O}_{{\mathbb P}(V)}(1)) = [E]$ is the hyperplane at
the infinity in ${\mathbb P}(V)$ due to the inclusion ${\mathcal
O}_\Gamma(-1)\oplus {\mathcal O}_\Gamma(-1) \longrightarrow {\mathcal
O}_\Gamma(-1)\oplus {\mathcal O}_\Gamma(-1)\oplus {\mathcal O}_\Gamma$. Therefore
we have
\begin{eqnarray*}
  c_1(X_2)\cdot \beta_2 = 3 |\eta|.
\end{eqnarray*}

By the definition of absolute Donaldson-Thomas invariants, we may
assume that
\begin{eqnarray*}
   c_1(X)\cdot \beta = \mbox{vdim} I_n(X,
   \beta) = \sum_{i=1}^r \deg ch_2(\gamma_{l_i}).
\end{eqnarray*}
Otherwise, the involved Donaldson-Thomas invariants of $X$ and
$\tilde{X}$ will vanish and the theorem holds.

We also have
\begin{eqnarray*}
c_1(X_1)\cdot
\beta_1=\text{vdim}I_n(X_1/E,\beta)=\sum_{i=1}^r\deg
ch_2(\gamma_{\ell_i})+\deg\epsilon_1^*\eta.
\end{eqnarray*}

 By Lemma
\ref{lem2.3}, we have
\begin{eqnarray*}
   c_1(X)\cdot \beta = c_1(X_1)\cdot \beta_1 +
   c_1(X_2)\cdot\beta_2 - 2 |\eta|.
\end{eqnarray*}
 Combining all the four equalities above,  we have
\begin{eqnarray*}
   0 = \deg \epsilon_1^*C_\eta + |\eta|
\end{eqnarray*}

Hence $|\eta| = 0$.

(i) Suppose that $\beta =m[\Gamma]$.  Notice that the virtual dimension of the
moduli space will be zero since $c_1(X)\cdot \beta=0$.  Let
$\Gamma_{\infty}$ be the curve coming from the inclusion $\mathcal
O_{\Gamma}\rightarrow V$, $F\cong\mathbb P^2$ be a fiber of $p$, $f$ be a
line in $F$. Then one can compute easily that
\begin{eqnarray*}
E\cdot \Gamma_{\infty}=0, \quad F\cdot \Gamma_{\infty}=1,\quad f\cdot F=0,
\quad f\cdot E=1.
\end{eqnarray*}
Therefore we can write $\beta_2=af+m[\Gamma_{\infty}]$. Since $E\cdot
\beta_2=0$, we have $a=0$. Therefore
$\beta_2=m [\Gamma_{\infty}]$ for some $m\ge 0$. Under the morphism
$\sigma$
\begin{eqnarray*}
\sigma\colon \mathcal X\longrightarrow \mathcal X\times \mathbb C\longrightarrow X,
\end{eqnarray*}
we have
$\beta=\sigma(\beta_1)+m[\Gamma]$ in $NE(X)$. $\beta_1$ can only be a union of
curves
$C_i$'s  not lying on $E$ and curves $D_j$'s on $E$. Since
$\mathbb R[\Gamma]$ is a ray, we must have $C_i=0$. For effective curves $D_j$ on
$E$, $D_j\cdot E\ne 0$. However since $\beta_1\cdot E=0$, we must have $D_j=0$.
Thus $\beta_1=0$.  Therefore, by the degeneration formula, we have
\begin{eqnarray}\label{Gamma}
Z'_{DT}(X;q)_{m[\Gamma]} =
Z'_{DT}({\tilde{X}}/{E};q)_0
\cdot Z'_{DT}({X_2}/{E};q)_{m[\Gamma_{\infty}]}=
 Z'_{DT}({X_2}/{E};q)_{m[\Gamma_{\infty}]}
\end{eqnarray}
\begin{eqnarray*}
Z'_{DT}(X^f;q)_{m[\Gamma^f]}  =
Z'_{DT}({\tilde{X^f}}/{E};q)_0
\cdot Z'_{DT}({X^f_2}/{E};q)_{m[\Gamma^f_{\infty}]}=
 Z'_{DT}({X^f_2}/{E};q)_{m[\Gamma^f_{\infty}]}.
\end{eqnarray*}

Observe that $(\tilde X_2, E)$ and $(\tilde X^f_2, E)$ are isomorphic. Therefore,
we  have
\begin{eqnarray*}
Z'_{DT}(X;q)_{m[\Gamma]}
 = Z'_{DT}(X^f;q)_{m[\Gamma_f]}.
\end{eqnarray*}

To write in another way for $\beta=m[\Gamma]$, we have
\begin{eqnarray*}
Z'_{DT}(X;q)_{\beta}
 = Z'_{DT}(X^f;q)_{-\phi_*(\beta)}.
\end{eqnarray*}

To prove (ii),   by the similar argument as in (i), we have
$\beta=\beta_1+m[\Gamma_{\infty}]$ with $m\ge 0$ and $\beta_1\cdot E=0$.

Furthermore, by the degeneration formula, we have

\begin{eqnarray}\label{flop-formula1}
  &&Z'_{DT}\big(X;q\mid
  \prod_{i=1}^r\tilde{\tau_0}(\phi^*\gamma_{l_i})\big)_\beta \nonumber\\
&=&
\sum_{ \beta=\beta_1+m[\Gamma_{\infty}],\atop
\beta_1\in \iota(H_2(X))}Z'_{DT}({\tilde{X}}/{E};q\mid
\prod_{i=1}^r\tilde{\tau_0}(\mu^*\phi^*\gamma_{l_i}))_{\beta_1}
\cdot Z'_{DT}({X_2}/{E};q)_{m[\Gamma_{\infty}]}.
\end{eqnarray}

Consider the map $c_*\colon H_2(X)=H_2(\mathcal X_t){\buildrel
i_{t*}\over \longrightarrow} H_2(\mathcal X){\buildrel i^{-1}_{0*}\over
\longrightarrow }H_2(\mathcal X_0)$. From Lemma 2.11 in \cite{L-R}, $c_*$ is
injective. Therefore we have
\begin{eqnarray*}
& &Z^{\prime}_{DT}\big(X; q, v \mid
  \prod_{i=1}^r\tilde{\tau_0}(\phi^*\gamma_{l_i})\big)
\\
&=&\sum_{\beta\in H_2(X)}Z'_{DT}\big(X;q\mid
   \prod_{i=1}^r\tilde{\tau_0}(\phi^*\gamma_{l_i})\big)_\beta v^{\beta}\\
&=&\sum_{\beta\in H_2(X)}\sum_{ \beta=\beta_1+m[\Gamma_{\infty}],\atop
\beta_1\in \iota(H_2(X))}Z'_{DT}({\tilde{X}}/{E};q\mid
\prod_{i=1}^r\tilde{\tau_0}(\mu^*\phi^*\gamma_{l_i}))_{\beta_1}v^{\beta_1}
\cdot Z'_{DT}({X_2}/{E};q)_{m[\Gamma_{\infty}]}v^{m[\Gamma_{\infty}]}\\
&=& \big(\sum_{\beta_1\in \iota(H_2(X))}Z'_{DT}({\tilde{X}}/{E};q\mid
\prod_{i=1}^r\tilde{\tau_0}(\mu^*\phi^*\gamma_{l_i}))_{\beta_1}v^{\beta_1}\big)
\cdot
\big(\sum_{m\ge
0}Z'_{DT}({X_2}/{E};q)_{m[\Gamma_{\infty}]}v^{m[\Gamma_{\infty}]}\big).
\end{eqnarray*}

Define a function $\Phi_X(q, v|\{\phi^*\gamma_{\ell_i}\})$ as follows
\begin{eqnarray*}
\Phi_X(q, v|\{\phi^*\gamma_{\ell_i}\})
=\sum_{\beta_1\in \iota(H_2(X))}Z'_{DT}({\tilde{X}}/{E};q\mid
\prod_{i=1}^r\tilde{\tau_0}(\mu^*\phi^*\gamma_{l_i}))_{\beta_1}v^{\beta_1}.
\end{eqnarray*}

Apply the formula (\ref{Gamma}) to $X=X_2$, we get
$Z'_{DT}({X_2}/{E};q)_{m[\Gamma_{\infty}]}=Z'_{DT}(X_2;q)_{m[\Gamma_{\infty}]} $.
We define a function $G(q, v, \Gamma_{\infty})$ as
follows:

\begin{eqnarray*}
&&G(q, v, \Gamma_{\infty})=\sum_{m\ge
0}Z'_{DT}({X_2}/{E};q)_{[m\Gamma_{\infty}]}v^{[m\Gamma_{\infty}]}
=\sum_{m\ge 0}Z'_{DT}(X_2;q)_{[m\Gamma_{\infty}]}v^{[m\Gamma_{\infty}]}\\
&=&Z'_{GW}(\mathcal O_{\Gamma_{\infty}}(-1)\oplus \mathcal
O_{\Gamma_{\infty}}(-1); u, v).
\end{eqnarray*}
The last equality is the Theorem 3 in \cite{MNOP1} for local Calabi-Yau $\mathcal O_{\Gamma_{\infty}}(-1)\oplus \mathcal
O_{\Gamma_{\infty}}(-1) $.

>From \cite{MNOP1}, we have
\begin{eqnarray*}
Z'_{GW}(\mathcal O_{\Gamma_{\infty}}(-1)\oplus \mathcal
O_{\Gamma_{\infty}}(-1) ; u, v)=exp\{F'_{GW}(\mathcal O_{\Gamma_{\infty}}(-1)\oplus \mathcal
O_{\Gamma_{\infty}}(-1) ; u, v)\},
\end{eqnarray*}

\begin{eqnarray*}
F'_{GW}=\sum_{d>0}\sum_{g\ge 0}N_{g, d}u^{2g-2}v^{d[\Gamma_{\infty}]},
\end{eqnarray*}
where $N_{g, d}$ is computed in \cite{F-P}:
\begin{eqnarray*}
N_{0, d}=\frac{1}{d^3},\quad N_{1, d}=\frac{1}{12d},\quad N_{g,
d}=\frac{|B_{2g}|d^{2g-3}}{2g\cdot (2g-2)!}\quad \hbox{for } g\ge 2.
\end{eqnarray*}
Therefore, we have
\begin{eqnarray*}
F'_{GW}=u^{-2}\sum_{d>0}\frac{1}{d^3}(v^{[\Gamma_{\infty}]})^d
+\sum_{d>0}\frac{1}{12d}(v^{[\Gamma_{\infty}]})^d
+\sum_{g\ge 2}\frac{|B_{2g}|}{2g\cdot
(2g-2)!}u^{2g-2}\sum_{d>0}d^{2g-3}(v^{[\Gamma_{\infty}]})^d.
\end{eqnarray*}

Now $
G(q, v,\Gamma_{\infty})/g(u,v,\Gamma_{\infty})$ has the analytic
continuation
\begin{eqnarray*}
exp\{\sum_{g\ge 2}\frac{|B_{2g}|}{2g\cdot
(2g-2)!}u^{2g-2} f_{2g-3}(v^{[\Gamma_{\infty}]})\}
\end{eqnarray*}
where $f_{2g-3}(x)$ is defined in the Lemma \ref{analytic-continuation}.

Applying the same argument to $X^f$, we also have
\begin{eqnarray}\label{flop-formula2}
  &&Z'_{DT}\big(X^f;q\mid
  \prod_{i=1}^r\tilde{\tau_0}(\gamma_{l_i})\big)_{\beta}\nonumber\\
&=&\sum_{ \beta=\beta_1+m[\Gamma^f_{\infty}],\atop
\beta_1\in \iota^f(H_2(X^f))}Z'_{DT}({\tilde{X}}/{E};q\mid
\prod_{i=1}^r\tilde{\tau_0}(\nu^*\gamma_{l_i}))_{\beta_1}
\cdot Z'_{DT}({X^f_2}/{E};q)_{m[\Gamma^f_{\infty}]}.
\end{eqnarray}
where $\nu\colon \widetilde X\rightarrow X^f$ is the blowup map,  $\tilde
X\cong\tilde{X^f}$.

Applying the same argument above  for $X$ to $X^f$,
define a function $\Phi_{X^f}(q, v|\{\gamma_{\ell_i}\})$ as follows
\begin{eqnarray*}
\Phi_X(q, v|\{\gamma_{\ell_i}\})
=\sum_{\beta_1\in \iota(H_2(X))}Z'_{DT}({\tilde{X}}/{E};q\mid
\prod_{i=1}^r\tilde{\tau_0}(\nu^*\gamma_{l_i}))_{\beta_1}v^{\beta_1}.
\end{eqnarray*}
We have (\ref{flop-split2}).

The function
$G(q, v,\Gamma^f_{\infty})/g(q, v,\Gamma^f_{\infty})$ has the analytic continuation
\begin{eqnarray*}
exp\{\sum_{g\ge 2}\frac{|B_{2g}|}{2g\cdot
(2g-2)!}u^{2g-2} f_{2g-3}(v^{[\Gamma^f_{\infty}]})\}
\end{eqnarray*}

>From  the Lemma \ref{analytic-continuation} and the fact that $\mu^*\phi^*=\nu^*$,
we proved (ii).
\end{proof}

 One should compare the Theorem \ref{flopthm} with Definition 1.1,
Theorem A and Corollary A.2 in \cite{L-R}. There, Li and Ruan
studied the question of naturality of quantum cohomology under birational
operations such as flops. They observed that one must use ananlytic continuation to
compare the quantum cohomology of two Calabi-Yau $3$-folds which are flop
equivalent. The similar phenomenon occurs for Donaldson-Thomas invariants. However,
there is a slight complexity due to the function $g(q, v,\Gamma)$ coming from genus
zero and genus one contributions. It is possible that genus zero and genus one
create an anomaly.

\end{document}